\documentclass[12pt, a4paper, oneside]{article}
\usepackage[top=1in, bottom=1.1in, left=0.9in, right=0.5in]{geometry}
\usepackage[parfill]{parskip}
\usepackage{amsmath}
\usepackage{amsthm}
\usepackage{graphicx}
\usepackage{amssymb}
\usepackage{epstopdf}
\usepackage{setspace}
\usepackage{authblk}
\usepackage{enumerate}
\usepackage{lmodern}
\usepackage[hypcap]{caption}
\usepackage[english]{babel}
\usepackage{fancyhdr}
\addto\captionsenglish{}
\addto\captionsenglish{}
\addto\captionsenglish{}
\addto\captionsenglish{}
\addto\captionsenglish{}

\newlength\longest

\newcommand{\R}{{\mathbb R}}

\begin{document}

\newtheorem{theorem}{Theorem}[section]
\newtheorem{lemma}[theorem]{Lemma}
\newtheorem{corollary}[theorem]{Corollary}
\newtheorem{proposition}[theorem]{Proposition}
\newtheorem{conjecture}[theorem]{Conjecture}
\newtheorem{problem}[theorem]{Problem}
\newtheorem{claim}[theorem]{Claim}
\theoremstyle{definition}
\newtheorem{assumption}[theorem]{Assumption}
\newtheorem{remark}[theorem]{Remark}
\newtheorem{definition}[theorem]{Definition}
\newtheorem{example}[theorem]{Example}
\theoremstyle{remark}
\newtheorem{notation}{Notasi}
\renewcommand{\thenotation}{}

\title{A Note on Inclusion Properties of Weighted Orlicz Spaces}
\author{Al Azhary Masta${}^{1}$, Ifronika${}^2$, Muhammad Taqiyuddin${}^3$}

\affil{${}^{1,3}$Department of Mathematics Education, Universitas Pendidikan Indonesia, Jl. Dr. Setiabudi 229, Bandung 40154}

\affil{${}^{2}$Analysis and Geometry Group, Faculty of Mathematics and Natural Sciences, Bandung Institute of Technology, Jl. Ganesha 10, Bandung 40132

${}^{1}$alazhari.masta@upi.edu, ${}^{2}$ifronika@math.itb.ac.id, ${}^{3}$taqi94@hotmail.com}
\date{}

\maketitle

\begin{abstract}
In this paper we present sufficient and necessary conditions for inclusion relation between two weighted Orlicz spaces which complete the Osan\c{c}liol result in 2014. One of the keys to prove our results is to use the norm of the characteristic functions of the balls in $\mathbb{R}^n$.
\bigskip

\end{abstract}

\section{Introduction}

Orlicz spaces are generalization of Lebesgue spaces which were firstly introduced by Z. W. Birnbaum and W. Orlicz in 1931 (see \cite{Luxemburg,Rao}). Let us first recall the definition of Orlicz spaces. Let $\Phi: [ 0, \infty ) \rightarrow [ 0, \infty ) $ be a Young function [that is,  $\Phi$ is convex, $\lim \limits_{t\to 0}\Phi(t)=0=\Phi(0)$, left-continuous and $ \lim \limits_{t\to\infty} \Phi(t) = \infty$], \textit{the Orlicz space} $L_{\Phi}(\mathbb{R}^n)$ is the set of measurable functions $f : \R^n \rightarrow \mathbb{R}$ such that $ \int_{\R^n} \Phi ( a|f(x)|) dx < \infty$
for some $a > 0$.  The space $L_\Phi(\R^n)$ is a Banach space equipped with the norm
$$
\| f \|_{L_\Phi(\R^n)} := \inf \left\{  {b>0:
	\int_{\R^n}\Phi \left(\frac{|f(x)|}{b} \right) dx \leq1}\right\}.
$$ 

Meanwhile, for $\Phi$ is a Young function, \textit{the weak Orlicz space} $wL_{\Phi}(\mathbb{R}^n)$ is the set of all measurable functions $f : \mathbb{R}^n \rightarrow \mathbb{R}$  such that
$$
\| f \|_{wL_{\Phi}(\mathbb{R}^n)} := \inf \left\{  {b>0:
	\mathop {\sup }\limits_{t > 0} \Phi(t) \mu\left( \Bigl\{ x \in \mathbb{R}^n : \frac{|f(x)|}{b} > t \Bigr\} \right)  \leq1}\right\} < \infty. $$

Now, we move to the weighted Orlicz spaces and weighted weak Orlicz spaces. Let $\Phi$ be a Young function and $u$ is a weight on $\R^n$ (i.e $u : \mathbb{R}^n \rightarrow ( 0, \infty )$ is a measurable function), \textit{the weighted Orlicz space} $L_{\Phi}^u(\mathbb{R}^n)$ is the set of all measurable functions $f : \mathbb{R}^n \rightarrow \mathbb{R}$  such that $uf \in L_{\Phi}(\mathbb{R}^n)$. Note that, the space $L_{\Phi}^u(\mathbb{R}^n)$ is a Banach space equipped with the norm $$\| f \|_{L_{\Phi}^u(\mathbb{R}^n)} := \| uf \|_{L_{\Phi}(\mathbb{R}^n)} = \inf \left\{  {b>0:
	\int_{\R^n}\Phi \left(\frac{|u(x)f(x)|}{b} \right) dx \leq1}\right\}.$$

Similar with weighted Orlicz spaces,  for a Young function $\Phi$ and a weight $u$ on $\R^n$, \textit{the weighted weak Orlicz space} $wL_{\Phi}^u(\mathbb{R}^n)$ is the set of all measurable functions $f : \mathbb{R}^n \rightarrow \mathbb{R}$  such that $\| f \|_{wL_{\Phi}^u(\mathbb{R}^n)} := \| uf \|_{wL_{\Phi}(\mathbb{R}^n)} < \infty $.

Let $u_1, u_2 : \mathbb{R}^n \rightarrow ( 0, \infty )$, we denote $u_1 \preceq u_2$ if there exists a constant $C > 0$ such that $u_1(x)  \leq C u_2(x)$ for all $ x \in \mathbb{R}^n$.  Note that, if $u_1 \preceq u_2$ then $\| f \|_{L_{\Phi}^{u_1}(\mathbb{R}^n)} \leq C \| f \|_{L_{\Phi}^{u_2}(\mathbb{R}^n)}$ and $\| f \|_{wL_{\Phi}^{u_1}(\mathbb{R}^n)} \leq C \| f \|_{wL_{\Phi}^{u_2}(\mathbb{R}^n)}.$

%Several authors have studied about Orlicz spaces, for examples see \cite{Lech,Luxemburg,Masta1,Orlicz}, etc. 

The study of Lebesgue spaces and Orlicz spaces has been studied by many researchers in the last few decades (see \cite{Kufner,Luxemburg,Lech,Masta1,Orlicz,Taqi}, etc.). In 1989, Maligranda \cite{Lech} discussed inclusion properties of Orlicz spaces. Later in 2016, Masta\textit{ et al.} \cite{Masta1} obtained sufficient and necessary conditions for inclusion relation between two Orlicz spaces and between two weak Orlicz spaces by using different technique from Maligranda. Moreover, they have found that two Orlicz spaces and two weak Orlicz spaces can be compared with respect to Young functions for any measurable set, although the Lebesgue space $L_p$ are not comparable with respect to the number $p$.

On the other hand, Osan\c{c}liol \cite{Alen} have proved sufficient and necessary conditions for inclusion relation between two weighted Orlicz spaces, as in the following theorem. 

\medskip

\begin{theorem} \cite{Alen}\label{Teorema1.2}
Let $\Phi$ be a continuous Young function satisfying the $\bigtriangleup_{2}$ condition [that is, there exists $ K > 0$ such that $\Phi(2t) \leq K \Phi(t)$ for all $ t \geq 0$], and $ u_1, u_2$ are measurable functions such that $u_i(x+y) \leq u_i(x) \cdot u_i(y)$ for every $ x, y \in \R^n$, where $i=1,2$. Then the following statements are equivalent:

{\parindent=0cm
	
	{\rm (1)}  $u_1 \preceq u_2$.
	
	{\rm (2)} $L_{\Phi}^{u_2}(\mathbb{R}^n) \subseteq L_{\Phi}^{u_1}(\mathbb{R}^n)$.
	
	{\rm (3)} There exists a constant $C>0$ such that $\|f\|_{L_{\Phi}^{u_1}(\mathbb{R}^n)} \leq C \|f\|_{L_{\Phi}^{u_2}(\mathbb{R}^n)}$, for every $ f \in L_{\Phi}^{u_2}(\mathbb{R}^n)$.
	
	\par}
\end{theorem}

Related result for weak type of Orlicz spaces can be found in \cite{Masta4}. 

In this paper, we are interested in studying the inclusion properties of weighted Orlicz spaces. In connection with Theorem \ref{Teorema1.2}, we shall prove inclusion relation between weighted Orlicz spaces with respect to Young functions $\Phi_1, \Phi_2$ and weights $u_1, u_2$.

To achieve our purpose, we will use the similar methods in \cite{Gunawan, Masta1, Masta2,Masta3,Alen} which pay attention to the characteristic functions of open balls in $\R^n$. Next, we recall some lemmas which will be used later in next section.

\medskip

\begin{lemma}\label{lemma:1.1} \cite{Nakai1}
Suppose that $\Phi$ is a Young function and $ \Phi^{-1}(s):=\inf \{r \geq 0 :
\Phi (r) > s \}$. We have

{\parindent=0cm
{\rm (1)} $\Phi^{-1}(0) = 0$.

{\rm (2)} $ \Phi^{-1}(s_1) \leq \Phi^{-1}(s_2)$ for  $s_1 \leq s_2$.

{\rm (3)} $\Phi (\Phi^{-1}(s)) \leq s \leq \Phi^{-1}(\Phi(s))$ for $0 \leq s <
\infty$.

\par}

\end{lemma}

\bigskip

\begin{lemma} \cite{Masta3}\label{lemma:1.2}
Let $\Phi_1, \Phi_2$ be Young functions. For any $s>0$, if there exists $C_1, C_2 > 0$ such that $\Phi^{-1}_2(s) \leq C_1 \Phi^{-1}_1(C_2s)$, then we have  $\Phi_1 (\frac{t}{C_1})  \leq C_2 \Phi_2(t)$ for $ t = \Phi^{-1}_2(s).$
\end{lemma}

In this paper, the letter $C$ will be used for constants that may change from line to line, while constants with subscripts, such as $C_{1}, C_{2}$, do not change in different lines.

\section{Results}

First,we will investigate the inclusion properties of weighted Orlicz spaces with respect to distinct Young functions $\Phi_1$ and $\Phi_2$. For getting the result, we give attention to estimate the norm of the characteristic function of open ball in $\R^n$ as in the following lemma.

\bigskip

\begin{lemma}\label{lemma:2.2}\cite{Ning,Masta1}
	Let $\Phi$ be a Young function, $a\in\R^n$, and $r >0$ be arbitrary. Then we have
	$ \left\|\frac{\chi_{B(a,r)}}{u} \right\|_{L_{\Phi}^u(\mathbb{R}^n)} = \|\chi_{B(a,r)}\|_{
		\L_\Phi(\mathbb{R}^n)} = \frac{1}{\Phi^{-1}\bigl( \frac{1}{|B(a,r)|} \bigr)}$ where $|B(a,r)|$ denotes the volume of open ball $B(a,r)$ centered at $a\in\R^n$ with radius $r >0$.
\end{lemma}

\medskip

Now we come to the inclusion relation between $L_{\Phi_1}^u(\mathbb{R}^n)$ and $L_{\Phi_2}^u(\mathbb{R}^n)$ with respect to Young functions $\Phi_1, \Phi_2$. Given two Young functions $\Phi_1, \Phi_2$, we write $\Phi_1 \prec \Phi_2$ if there exists a constant $C > 0$ such that $\Phi_1(t) \leq \Phi_2(Ct)$ for all $ t > 0$.    

\medskip
  
\begin{theorem}\label{theorem:2.3}
	Let $\Phi_1, \Phi_2$ be Young functions and $ u: \mathbb{R}^n \rightarrow (0,\infty)$ be a measurable function. Then the following statements are equivalent:
	
	{\parindent=0cm
		{\rm (1)} $\Phi_1 \prec \Phi_2$.
		
		{\rm (2)} $L_{\Phi_2}^u(\mathbb{R}^n) \subseteq L_{\Phi_1}^u(\mathbb{R}^n)$.
		
		{\rm (3)} There exists a constant $C > 0$ such that $\| f \|_{L_{\Phi_1}^u(\mathbb{R}^n)} \leq C\| f \|_{L_{\Phi_2}^u(\mathbb{R}^n)}$, for every $ f \in L_{\Phi_2}^u(\mathbb{R}^n)$.
		
		\par}
\end{theorem}

\noindent{\it Proof}.
Assume that (1) holds. Suppose that  $ f \in L_{\Phi_2}^u(\mathbb{R}^n)$.
Observe that
\[
\int_{\mathbb{R}^n}\Phi_1 \left( \frac{|u(x)f(x)|}{C \| f \|_{L_{\Phi_2}^u(\mathbb{R}^n)}} \right) dx
\leq \int_{\mathbb{R}^n}\Phi_2 \left(\frac{C|u(x)f(x)|}{ C\| f \|_{L_{\Phi_2}^u(\mathbb{R}^n)}} \right)dx=
\int_{\mathbb{R}^n} \Phi_2 \left( \frac{|u(x)f(x)|}{ \| f \|_{L_{\Phi_2}^u(\mathbb{R}^n)}} \right) dx\le 1.
\]
By definition of $\| \cdot \|_{L_{\Phi}^u(\mathbb{R}^n)}$, we have $\| f \|_{L_{\Phi_1}^u(\mathbb{R}^n)}
\leq C \| f \|_{L_{\Phi_2}^u(\mathbb{R}^n)}.$
This proves that $L_{\Phi_2}^u(\mathbb{R}^n) \subseteq L_{\Phi_1}^u(\mathbb{R}^n)$.

Next, since $(L_{\Phi_1}^u(\mathbb{R}^n), L_{\Phi_2}^u(\mathbb{R}^n))$ is a Banach pair,
it follows from \cite[Lemma 3.3]{Krein} that (2) and (3) are equivalent. It thus remains to show that
(3) implies (1).

Assume now that (3) holds. By Lemma \ref{lemma:2.2}, we have
$$
\frac{1}{\Phi_1^{-1}\Bigl(\frac{1}{|B(a,r)|}\Bigr)} = \left\|  \frac{\chi_{B(a,r)}}{u} \right\|_{L_{\Phi_1}^u(\mathbb{R}^n)}
\leq C \left\|\frac{\chi_{B(a,r)}}{u} \right\|_{L_{\Phi_2}^u(\mathbb{R}^n)} =  \frac{C}{\Phi_2^{-1}\Bigl(\frac{1}{|B(a,r)|}\Bigr)}. 
$$
 
Since $\frac{1}{\Phi_1^{-1}\Bigl(\frac{1}{|B(a,r)|}\Bigr)}  \leq \frac{C}{\Phi_2^{-1}\Bigl(\frac{1}{|B(a,r)|}\Bigr)}$ is equivalent to $C \Phi_1^{-1}(\frac{1}{|B(a,r)|}) \geq \Phi_2^{-1}(\frac{1}{|B(a,r)|})$
for arbitrary $a\in\mathbb{R}^n$ and $r > 0$, by Lemma \ref{lemma:1.2}, we have
$$
\Phi_1\Bigl(\frac{t_0}{C}\Bigr) \leq  \Phi_2(t_0),
$$
for $t_0 = \Phi_2^{-1}( \frac{1}{|B(a,r)|})$. Since $a\in\mathbb{R}^n$ and $r > 0$ are arbitrary, we conclude that $\Phi_1(t) \leq \Phi_2(Ct)$ for every $t > 0$.\qed 

\medskip

\remark For $u(x)=1$, Theorem \ref{theorem:2.3} reduces to Theorem 2.5 in \cite{Masta1}. 

\medskip
Next, we also give the sufficient and necessary conditions for inclusion relation between weighted Orlicz spaces $L_{\Phi_1}^{u_1}(\mathbb{R}^n)$ and $L_{\Phi_2}^{u_2}(\mathbb{R}^n)$ with respect to Young functions $\Phi_1, \Phi_2$ and weights $u_1, u_2$. To get the result, we need the following lemma.

\medskip

\begin{lemma}\label{lemma:2.6}
Let $u : \mathbb{R}^n \rightarrow (0,\infty)$ be a measurable function such that $u(x+y) \leq u(x) \cdot u(y)$ for every $ x,y \in \R^n$. If $\Phi$ is a Young function, then:
{\parindent=0cm

{\rm (1)}  For all $ f \in L_{\Phi}^{u}(\mathbb{R}^n)$ and for all $ x \in \mathbb{R}^n$, we have $\|T_{x}f\|_{L_{\Phi}^{u}(\mathbb{R}^n)} \leq u(x)\|f\|_{L_{\Phi}^{u}(\mathbb{R}^n)},$ where $T_{x}f(y)=f(y-x)$.

{\rm (2)} If $ f \in L_{\Phi}^{u}(\mathbb{R}^n)$ and $f \neq 0$, then there exists a constant $ C > 0$ (depends on $f$) such that $$\frac{u(x)}{C} \leq \|T_{x}f\|_{L_{\Phi}^{u}(\mathbb{R}^n)} \leq C u(x).$$

\par}
\end{lemma}
\noindent{\it Proof}.

{\parindent=0cm
	
	{\rm (1)} Let $ f \in L_{\Phi}^{u}(\mathbb{R}^n)$ and $T_{x}f(y)=f(y-x)$, then $$\int_{\R^n}\Phi \Bigl(\frac{|u(v)f(v)|}{\|f\|_{L_{\Phi}^{u}(\mathbb{R}^n)}} \Bigr) dv \leq 1.$$ Observe that (by setting $v:=y-x$), we have
	\begin{align*}
 \int_{\R^n}\Phi \Bigl(\frac{|u(y)T_{x}f(y)|}{u(x)\|f\|_{L_{\Phi}^{u}(\mathbb{R}^n)}} \Bigr) dy &= \int_{\R^n}\Phi \Bigl(\frac{|u(y)f(y-x)|}{u(x)\|f\|_{L_{\Phi}^{u}(\mathbb{R}^n)}} \Bigr) dy\\
	& =  \int_{\R^n}\Phi \Bigl(\frac{|u(v+x)f(v)|}{u(x)\|f\|_{L_{\Phi}^{u}(\mathbb{R}^n)}} \Bigr) dv\\
	& \leq  \int_{\R^n}\Phi \Bigl(\frac{|u(v)u(x)f(v)|}{u(x)\|f\|_{L_{\Phi}^{u}(\mathbb{R}^n)}} \Bigr) dv\\
	& =  \int_{\R^n}\Phi \Bigl(\frac{|u(v)f(v)|}{\|f\|_{L_{\Phi}^{u}(\mathbb{R}^n)}} \Bigr) dv\\
	& \leq 1.
	\end{align*}
	
	This shows that $\|T_{x}f\|_{L_{\Phi}^{u}(\mathbb{R}^n)}\leq u(x)\|f\|_{L_{\Phi}^{u}(\mathbb{R}^n)}$.
	
	{\rm (2)} Let $ f \in L_{\Phi}^{u}(\mathbb{R}^n)$ and $f\neq0$, then there exist a constant $ C > 0$ (depends on $f$) such that $\|f\|_{L_{\Phi}^{u}(\mathbb{R}^n)}\leq C.$ By Lemma \ref{lemma:2.6} (1), then we have $$\|T_{x}f\|_{L_{\phi,\Phi}^{u}(\mathbb{R}^n)} \leq C u(x),$$ for every $ x \in \R^n$.
	
Since $f(x) \neq 0$ for every $ x\in \R^n$,we have $\|T_{x}f\|_{L_{\Phi}^{u}(\mathbb{R}^n)} > 0$. Observe that\\
	\begin{align*}
	\int_{\R^n}\Phi \Bigl(\frac{|u(x)f(v)|}{\mathop {\sup }\limits_{ v \in \R^n}u(-v)\|T_{x}f\|_{L_{\Phi}^{u}(\mathbb{R}^n)}} \Bigr) dv	& \leq \int_{\R^n}\Phi \Bigl(\frac{|u(x)f(v)|}{u(-v)\|T_{x}f\|_{L_{\Phi}^{u}(\mathbb{R}^n)}} \Bigr) dv\\
	& \leq \int_{\R^n}\Phi \Bigl(\frac{|u(v+x)f(v)|}{\|T_{x}f\|_{L_{\Phi}^{u}(\mathbb{R}^n)}} \Bigr) dv\\
	& \leq \int_{\R^n}\Phi \Bigl(\frac{|u(y)f(y-x)|}{\|T_{x}f\|_{L_{\Phi}^{u}(\mathbb{R}^n)}} \Bigr) dy\\
	& = \int_{\R^n}\Phi \Bigl(\frac{|u(y)L_{x}f(y)|}{\|T_{x}f\|_{L_{\Phi}^{u}(\mathbb{R}^n)}} \Bigr) dy\\
	& \leq 1.
	\end{align*}
	
	This shows that $\frac{u(x)\|f\|_{_{L_{\Phi}(\mathbb{R}^n)}}}{\mathop {\sup }\limits_{ v \in R^n}u(-v)} \leq \|T_{x}f\|_{L_{\Phi}^{u}(\mathbb{R}^n)}$.
	
	Choose, $C:=\mathop {\max } \Bigl\{C_1, \frac{\mathop {\sup }\limits_{ v \in \R^n }u(-v)}{\|f\|_{L_{\Phi}(\mathbb{R}^n)}} \Bigr\}$. Hence, we conclude that $\frac{u(x)}{C} \leq \|T_{x}f\|_{L_{\Phi}^{u}(\mathbb{R}^n)}\leq Cu(x)$, as desired. \qed
	\par}

\medskip

Now, we present the sufficient and necessary conditions for the inclusion properties of weighted Orlicz spaces $L_{\Phi_1}^{u_1}(\mathbb{R}^n)$ and $L_{\Phi_2}^{u_2}(\mathbb{R}^n)$ with respect to Young functions $\Phi_1, \Phi_2$ and weights $u_1, u_2$.

\medskip

\begin{theorem}\label{theorem:2.4}
Let $\Phi_{1}, \Phi_{2}$ be Young functions such that $\Phi_1 \prec \Phi_2$ and $ u_1, u_2$ are measurable functions such that $u_i(x+y) \leq u_i(x) \cdot u_i(y)$ for every $ x, y \in \R^n$, where $i=1,2$. Then the following statements are equivalent:
	
	{\parindent=0cm
		
		{\rm (1)}  $u_1 \preceq u_2$.
		
		{\rm (2)} $L_{\Phi_2}^{u_2}(\mathbb{R}^n) \subseteq L_{\Phi_1}^{u_1}(\mathbb{R}^n)$.
		
		{\rm (3)} There exists a constant $C>0$ such that $\|f\|_{L_{\Phi_1}^{u_1}(\mathbb{R}^n)} \leq C \|f\|_{L_{\Phi_2}^{u_2}(\mathbb{R}^n)}$, for every $ f \in L_{\Phi_2}^{u_2}(\mathbb{R}^n)$.
		
		\par}
\end{theorem}

Proof.

Assume that (1) holds. Let  $ f$ be an element of $ L_{\Phi_2}^{u_2}(\mathbb{R}^n)$. Since $\Phi_1 \prec \Phi_2$ and $u_1 \preceq u_2$, there exists constants $C_1, C_2>0$ such that $\Phi_1(t) \leq \Phi_2(C_1t)$ for all $t > 0$ and $ u_1(x) \leq C_2 u_2(x)$ for every $ x \in \R^n$. Using a similar argument in the proof of Theorem \ref{theorem:2.3} we have $$\|  f \|_{L_{\Phi_1}^{u_1}(\R^n)} \leq C_1\|  f \|_{L_{\Phi_2}^{u_1}(\R^n)} \leq C_1C_2 \|  f \|_{L_{\Phi_2}^{u_2}(\R^n)}.$$ 

As before, we have that (2) and (3) are equivalent. It thus remains to show that (3) implies (1).
Assume that (3) holds. By Lemma \ref{lemma:2.6}, we have

$$\frac{u_{1}(x)}{C}\leq \|T_{x}f\|_{L_{\Phi_1}^{u_1}(\mathbb{R}^n)} \leq C \|T_{x}f\|_{L_{\Phi_2}^{u_2}(\mathbb{R}^n)}\leq Cu_2(x),$$ for every $x \in \mathbb{R}^n$. So, we obtain $u_1 \preceq u_2$.\qed

\bigskip
Note that, for $\Phi_1(x) = \Phi_2(x)$ for every $ x \in \R^n$, Theorem \ref{theorem:2.4} reduces to Theorem \ref{Teorema1.2}.

\medskip

\remark It follows from Theorems \ref{theorem:2.3} and \ref{theorem:2.4} that there cannot be an inclusion relation between $L_{p_1}^{u_1}(\R^n)$ and $L_{p_2}^{u_2}(\R^n)$ for distinct values of $p_1$ and $p_2$. In spite of that, for finite measure set $X$ we can obtain inclusion relation between $L_{p_1}^{u_1}(X)$ and $L_{p_2}^{u_2}(X)$ presented in the next section.

\section{An additional case}

In the following, we will give sufficient condition for H\"{o}lder's inequality in weighted Orlicz spaces which will be used to obtain inclusion relation between $L_{p_1}^{u_1}(X)$ and $L_{p_2}^{u_2}(X)$.   

\medskip

\begin{theorem}\label{theorem:3.1} (H\"{o}lder's inequality)
Let $X$ be a measurable set, $\Phi_1, \Phi_2$, $\Phi_3$ be Young functions and $u_1, u_2, u_3: X \rightarrow \R$ be measurable functions such that $\Phi^{-1}_1(t)\Phi^{-1}_2(t) \leq \Phi^{-1}_3(t)$ for every  $t > 0$
and  $u_3(x) \leq u_1(x)u_2(x)$ for every $ x \in X$. If $f_1 \in L_{\Phi_1}^{u_1}(X)$ and $f_2 \in L_{\Phi_2}^{u_2}(X)$, then $f_1 f_2 \in L_{\Phi_3}^{u_3}(X)$ with $$\| f_1 f_2 \|_{L_{\Phi_3}^{u_3}(X)} \leq 2 \| f_1 \|_{L_{\Phi_1}^{u_1}(X)} \| f_2 \|_{L_{\Phi_2}^{u_2}(X)}.$$
\end{theorem}

\noindent{\it Proof}. 

Let $s,t \ge 0$.
Without loss of generality, suppose that $\Phi_1(s) \leq \Phi_2(t)$. By Lemma \ref{lemma:1.1}(3), we obtain
$$
st \leq \Phi^{-1}_1(\Phi_1(s))\Phi^{-1}_2(\Phi_2(t))\leq \Phi^{-1}_1(\Phi_2(t))
\Phi^{-1}_2(\Phi_2(t))\leq \Phi^{-1}_3(\Phi_2(t)).
$$
Hence $\Phi_3(st)\leq \Phi_3(\Phi^{-1}_3(\Phi_2(t))) \leq \Phi_2(t) \leq
\Phi_2(t) + \Phi_1(s)$. Since $\Phi$ is a convex function, we have

\begin{align*}
\int_{X}\Phi_3 \left( \frac{|u_3(x)f_1(x)f_2(x)|}{2 \| f_1 \|_{L_{\Phi_1}^{u_1}(X)}
	\| f_2 \|_{L_{\Phi_2}^{u_2}(X)}} \right) dx &\leq \frac{1}{2} \int_{X}\Phi_3 \left(
\frac{|u_3(x)f_1(x)f_2(x)|}{\| f_1 \|_{L_{\Phi_1}^{u_1}(X)} \| f_2 \|_{L_{\Phi_2}^{u_2}(X)}}\right) dx \\
& \leq \frac{1}{2} \int_{X}\Phi_3 \left(
\frac{|u_1(x)u_2(x)f_1(x)f_2(x)|}{\| f_1 \|_{L_{\Phi_1}^{u_1}(X)} \| f_2 \|_{L_{\Phi_2}^{u_2}(X)}}\right) dx.
\end{align*}

On the other hand, by Lemma \ref{lemma:2.2} we obtain\\
$$ \int_{X}\Phi_3 \left(
\frac{|u_1(x)u_2(x)f_1(x)f_2(x)|}{\| f_1 \|_{L_{\Phi_1}^{u_1}(X)} \| f_2 \|_{L_{\Phi_2}^{u_2}(X)}}\right) dx \leq \int_{X}\Phi_1 \left( \frac{|u_1(x)f_1(x)|}{\| f_1 \|_{L_{\Phi_1}^{u_1}(X)}}
\right) dx + \int_{X}\Phi_2 \left( \frac{|f_2(x)|}{\| f_2 \|_{L_{\Phi_2}^{u_2}(X)}}
\right) dx \leq 2,$$

whenever $f_1 \in L_{\Phi_1}^{u_1}(X)$ and $f_2 \in L_{\Phi_2}^{u_2}(X)$. By the
definition of $\| \cdot \|_{L_{\Phi_3}^{u_1}(X)}$, we have $\| f_1f_2 \|_{L_{\Phi_3}^{u_3}(X)} \leq
2 \| f_1 \|_{L_{\Phi_1}^{u_1}(X)}\|  f_2 \|_{L_{\Phi_2}^{u_2}(X)}$, as desired.\qed

\bigskip

\begin{corollary}\label{corollary:3.2}
	Let $X:=B(a,r_0) \subseteq \mathbb{R}^n$ for some $a\in\mathbb{R}^n$ and $r_{0} > 0$.
	If $\Phi_1, \Phi_2$ are
	two Young functions, $u_1, u_2: X \rightarrow \R$ are measurable functions and there are a Young function $\Phi$ and a weight $0 < u(x) \leq 1$ for every $ x \in X$ such that
	$
	\Phi^{-1}_1(t)\Phi^{-1}(t) \leq \Phi^{-1}_2(t)
	$
	for every $ t \geq 0$ and $ u_1(x) \leq u(x)u_2(x)$ for every $ x\in X$, then $$L_{\Phi_1}^{u_1}(X)\subseteq L_{\Phi_2}^{u_2}(X)$$ with
	$
	\| f \|_{L_{\Phi_2}^{u_2}(X)} \leq \frac{2}{\Phi^{-1}(\frac{1}
		{|B(a,r_{0})|})} \| f \|_{L_{\Phi_1} ^{u_1}(X)}
	$
	for $ f \in L_{\Phi_1}^{u_1}(X)$.
\end{corollary}

\noindent{\it Proof}. Since $ 0 < u(x) \leq 1$ for every $ x \in X$, we have $ \|f\|_{L_{\Phi_2}^{u_2}(X)} \leq \|\frac{f}{u}\|_{L_{\Phi_2}^{u_2}(X)}$. 
Let $f \in L_{\Phi_1}^{u_1}(X)$, by Theorem \ref{theorem:3.1} and choosing $g:=\chi_{B(a,r_0)}$,
we obtain
\begin{align*}
\| f\|_{L_{\Phi_2}^{u_2}(X)} = & \| f \chi_{B(a,r_0)} \|_{L_{\Phi_2}^{u_2}(X)} \\
= & \| f g \|_{L_{\Phi_2}^{u_2}(X)} \\
\leq & \left\|\frac{fg}{u} \right\|_{L_{\Phi_2}^{u_2}(X)} \\  
\leq & 2 \left\|\frac{g}{u}
\right\|_{L_{\Phi}^u(X)} \| f \|_{L_{\Phi_1}^{u_1}(X)} \\
= &\frac{2}{\Phi^{-1}(\frac{1}{|B(a,r_{0})|})} \| f \|_{L_{\Phi_1}^{u_1}(X)}.
\end{align*}

This shows that $L_{\Phi_1}^{u_1}(X) \subseteq L_{\Phi_2}^{u_2}(X)$. \qed

\medskip

We shall now discuss the inclusion properties of weighted weak Lebesgue spaces $L_{p_1}^{u_1}(X)$ and $L_{p_2}^{u_2}(X)$ with respect to distinct values of $p_1$ and $p_2$ as well as $u_1$ and $u_2$.

\medskip

\begin{corollary}
	Let $X:=B(a,r_0)$ for some $a\in\mathbb{R}^n$ and $r_0>0$. If $ 1 \leq p_{2} <
	p_{1} < \infty$ and $u_1, u_2: X \rightarrow \R$ are measurable functions such that $u_1(x) \leq u_2(x)$ for every $ x \in X$, then $$L_{p_{1}}^{u_1}(X) \subseteq L_{p_{2}}^{u_2}(X).$$
\end{corollary}

\noindent{\it Proof}.
Let $\Phi_1(t):= t^{p_1}, \Phi_2(t):= t^{p_2}$, $\Phi(t):=
t^{\frac{p_1 p_2}{p_1- p_2}}$ for every $t\ge 0$. Since $ 1 \leq p_{2} < p_{1} < \infty$,
we have $\frac{p_1 p_2}{p_1- p_2} > 1$. Thus, $\Phi_1,\
\Phi_2$, and $\Phi$ are Young functions. Now, define $u(x) = \frac{u_1(x)}{u_2(x)}$ for every $x \in X$. Observe that, using the definition
of $\Phi^{-1}$ and Lemma \ref{lemma:1.1}, we have
$$
\Phi^{-1}_{1}(t) =t^{\frac{1}{p_1}},\ \Phi^{-1}_{2}(t) =t^{\frac{1}{p_2}},\
{\rm and}\ \Phi^{-1}(t) =t^{\frac{p_1- p_2}{p_1p_2}}.
$$
Moreover, $\Phi^{-1}_{1}(t)\Phi^{-1}(t) = t^{\frac{1}{p_1}}t^{\frac{p_1-
		p_2}{p_1p_2}} = t^{\frac{1}{p_2}} = \Phi^{-1}_{2}(t)$ and $ u_1(x) = \frac{u_1(x)}{u_2(x)}u_2(x) = u(x)u_2(x)$. So it follows
from Corollary \ref{corollary:3.2} that $\| f \|_{L_{p_2}^{u_2}(X)} \leq 2|B(a,r_0)|^{\frac{p_1 - p_2}{p_1p_2}}  \| f \|_{L_{p_1}^{u_1}(X)}$, and therefore we can conclude that $L_{p_{1}}^{u_1}(X)
\subseteq L_{p_{2}}^{u_2}(X)$. \qed

\section{Concluding Remarks}

We have shown the inclusion properties of weighted Orlicz spaces for distinct Young functions $\Phi_1, \Phi_2$ and weights $u_1, u_2$.  The inclusion properties of weighted Orlicz spaces are generalization of inclusion properties of Orlicz spaces in \cite{Masta1} and inclusion properties of weighted Lebesgue spaces. In the proof of our results, we used the norm of characteristic function on $\R^n$ and estimated the norm of the translation functions in $\R^n$. 

Furthemore, from Theorem \ref{theorem:2.4} and Lemma 1.1, Theorem 2.8 in \cite{Masta4}, we also have the following inclusion
relations
\[
\begin{array}{ccc}
L^{u_2}_{\Phi_2} & \rightarrow &L^{u_1}_{\Phi_1}\\
\downarrow & \searrow & \downarrow\\
wL^{u_2}_{\Phi_2} &\rightarrow & wL^{u_1}_{\Phi_1}
\end{array}
\]
for $\Phi_1 \prec \Phi_2$ and $u_1 \preceq u_2$, where the arrows
mean `contained in' or `embedded into'.

\medskip

\textbf{Acknowledgement}. The first author is supported by Hibah Penguatan Kompetensi UPI 2018. The authors thank the referee for his/her useful remarks on the earlier version of this paper.


\begin{thebibliography}{10}

%\bibitem{H.G. Feichtinger}
%H.G. Feichtinger, ``Gewichtsfunktionen auf lokalkompakten gruppen'', \emph{\"{O}stereich. Akad. Wiss. Math. - Natur. Kl. Sitzungsber. II} \textbf{188}8-10 (1979), 451--471. 
%


\bibitem {Gunawan}
H. Gunawan, D.I. Hakim, K.M. Limanta, and A.A. Masta, ``Inclusion properties of generalized Morrey spaces'', \emph{Math. Nachr.} \textbf{290} (2017), 332--340. [DOI: 10.1002/mana.201500425]

\bibitem {Krein}
S.G Kre\v{i}n, Yu.\={I} Petun\={i}n, and E.M. Sem\"{e}nov, \emph{Interpolation of Linear
Operators}, Translation of Mathematical Monograph vol. 54, American Mathematical Society, Providence, R.I., 1982.


\bibitem{Kufner}
A. Kufner, O. John, and S. Fu\"{c}ik, \emph{Function Spaces}, Noordhoff International
Publishing, Czechoslovakia, 1977.

\bibitem{Ning}
N. Liu and Y. Ye, ``Weak Orlicz space and its convergence theorems'', \emph{Acta Math.
	Sci. Ser. B} {\bf 30}-5 (2010), 1492--1500.

\bibitem{Luxemburg}
W.A.J. Luxemburg, \emph{Banach Function Spaces}, Thesis, Technische Hogeschool te Delft, 1955.

\bibitem{Lech} L. Maligranda, \emph{Orlicz Spaces and Interpolation}, Departamento de Matem\'{a}tica, Universidade Estadual de Campinas, 1989.

\bibitem{Masta1} A.A. Masta, H. Gunawan, and W. Setya-Budhi, ``Inclusion property of Orlicz and weak Orlicz spaces'', \emph{J. Math. Fund. Sci}. \textbf{48}-3 (2016), 193--203 [DOI: http://dx.doi.org/10.56142Fj.math.fund.sci.2016.48.3.1].

\bibitem{Masta2} A.A. Masta, H. Gunawan, and W. Setya-Budhi, ``An inclusion property of Orlicz-Morrey spaces'',  \emph{J. Phys.: Conf. Ser.}, \textbf{893} 012015 (2017), 1--8 [DOI: https://doi.org/10.1088/1742-6596/893/1/012015].

\bibitem{Masta3} A.A. Masta, H. Gunawan, and W. Setya-Budhi, ``On Inclusion Properties of Two Versions of Orlicz-Morrey Spaces'',  \emph{Mediterr. J. Math.}, \textbf{14}-6 (2017), 228--239. [DOI: https://doi.org/10.1007/s00009-017-1030-7]

\bibitem{Masta4} A.A. Masta, Ifronika, and M. Taqiyuddin, ``Inclusion properties of weighted weak Orlicz spaces'',  research report (https://arxiv.org/abs/1710.04537), 2017.


\bibitem{Nakai1}	
E. Nakai, ``On Orlicz-Morrey spaces'', research report
[http://repository.kulib.kyoto-u.ac.jp/dspace/bitstream/2433/58769/1/1520-
10.pdf,
accessed on August 17, 2015.]

\bibitem{Orlicz}
W. Orlicz, \emph{Linear Functional Analysis (Series in Real Analysis Volume 4)},
World Scientific, Singapore, 1992.

\bibitem{Alen}
A. Osan\c{c}liol, ``Inclusion between weighted Orlicz spaces'', \emph{J. Inequal. Appl.}
\textbf{2014}-390 (2014), 1--8 [DOI: https://doi.org/10.1186/1029-242X-2014-390].

%\bibitem{Oztop}
%S. \"{O}ztop and A.T. G\"{u}rkanli, ``Multipliers and tensor products of weighted $L^p$ spaces'', \emph{Acta Math. Sci.}
%\textbf{201}-1 (2001), 41--49.	.

\bibitem{Rao}
M.M. Rao and Z.D. Ren, \emph{Theory of Orlicz spaces, volume 146 of Monographs and
	Textbooks in Pure and Applied Mathematics}, Marcel Dekker, Inc., New York, 1991.

\bibitem{Taqi}
M. Taqiyuddin and A.A. Masta, ``Inclusion properties of Orlicz spaces and weak Orlicz spaces generated by concave function'', \textit{IOP Conf. Ser.: Mater. Sci. Eng.,} \textbf{288} 012103 (2018), 1--5 [DOI: https://doi.org/10.1088/1757-899X/288/1/012103].

\bibitem{Xueying}
X. Zhang and C. Zhang, ``Weak Orlicz spaces generated by concave functions'',
International Conference on Information Science and Technology (ICIST) 2011, 42--44.



\end{thebibliography}
\end{document}